\documentclass[a4paper,12pt]{article}
\usepackage[T2A]{fontenc} 
\usepackage[utf8]{inputenc}
\usepackage[english]{babel}

\usepackage{mathrsfs}
\usepackage[top=2cm,left=3cm,right=1.5cm,bottom=2cm]{geometry}
\usepackage{amsmath}
\usepackage{amsfonts}
\usepackage{amssymb}
\usepackage{amsthm} 
\usepackage{enumerate}
\usepackage{hyperref}
\usepackage{alphalph}
\usepackage[subnum]{cases}
\usepackage{array}

\numberwithin{equation}{subsection}

\theoremstyle{plain}
\newtheorem{theorem}{Theorem}[section]

\theoremstyle{definition}

\theoremstyle{remark}
\newtheorem{remark}[theorem]{Remark}

\numberwithin{equation}{section}

\begin{document}


\begin{center}\textbf{About one inverse problem for a Hill's equation with double eigenvalues}\end{center}
\vskip 0.3 cm
\begin{center}\textbf{B.\,N. Biyarov}\end{center}
\\

%
\begin{abstract} 
In this paper, we establish necessary and sufficient conditions for the doubleness all of the eigenvalues, except the lowest, periodic and anti-periodic problems for Hill's equation in terms of the complex-valued potential $q(x)$.
\end{abstract}

\section{Introduction}
\label{sec:1} 
In the present paper we study the Hill's operator
\[\widehat{L}_{H} = -\dfrac{d^2}{dx^2} + \ q(x), \]
in the Hilbert space $L^2(0, 1)$, where $\ q(x)=q(x+1)$ is an arbitrary complex-valued function of $L^2(0, 1)$.
The closure in $L^2(0, 1)$ of the operator $\widehat{L}_{H}$ considered on $C^{\infty}[0, 1]$ is the maximal operator $\widehat{L}_{H}$ with the domain
\[D(\widehat{L}_{H}) = \{y \in L^2(0, 1): \ y,\ y' \in AC [0, 1], \  y'' - q(x)y \in L^2(0, 1)\}. \]
We consider the operator $L_D=\widehat{L}_{H}$ on the domain
\[D(L_D) = \{y \in D(\widehat{L}_{H}): \ y(0)=y(1)=0 \}, \]
the operator $L_N=\widehat{L}_{H}$ on the domain
\[D(L_N) = \{y \in D(\widehat{L}_{H}): \ y'(0)=y'(1)=0 \}, \]
the operator $L_{DN}=\widehat{L}_{H}$ on the domain
\[D(L_{DN}) = \{y \in D(\widehat{L}_{H}): \ y(0)=y'(1)=0 \}, \]
and the operator $L_{ND}=\widehat{L}_{H}$ on the domain
\[D(L_{ND}) = \{y \in D(\widehat{L}_{H}): \ y'(0)=y(1)=0 \}. \]
Also we consider the operator $L_P = \widehat{L}_{H}$ on the domain
\[
D(L_P)=\{ y \in D(\widehat{L}_{H}): y(0)=y(1), y'(0)=y'(1) \},
\]
the operator $L_{AP}=\widehat{L}_{H}$ on the domain
\[
D(L_{AP})=\{ y \in D(\widehat{L}_{H}): y(0)=-y(1), y'(0)=-y'(1) \},
\]
the operator $L_{D(\frac{1}{2})}=\widehat{L}_H$ on $[0,\frac{1}{2}]$ on the domain
\[
    D(L_{D(\frac{1}{2})}) = \{ y \in D(\widehat{L}_H): y(0)=y(1/2)=0 \},
\]
and  the operator $L_{N(\frac{1}{2})}=\widehat{L}_H$ on $[0, \frac{1}{2}]$ on the domain
\[
    D(L_{N(\frac{1}{2})})=\{ y \in  D(\widehat{L}_H): y'(0)=y'(1/2)=0 \}. 
\] 
Here we use subscripts D, N, DN, ND, P and AP meaning Dirichlet, Neumann, Dirichlet-Neumann, Neumann-Dirichlet, Periodic and Anti-Periodic operators, respectivelly. By $\sigma(A)$ we denote the spectrum of the operator $A$.

Consider the Hill's equation in the Hilbert space $L^2(0, 1)$
\begin{equation}\label{eq:1.1}
    \widehat{L}_Hy \equiv - y'' + q(x)y = \lambda^2y,
\end{equation}
where $\ q(x)=q(x+1)$ is the complex-valued function of $L^2(0, 1)$.
By $c(x, \lambda)$ and $s(x, \lambda)$ we denote the fundamental system of solutions to the equation \eqref{eq:1.1} corresponding to the initial conditions 
\[c(0, \lambda) = s'(0, \lambda)=1 \;\, \mbox{and} \;\,  c'(0, \lambda) = s(0, \lambda)=0.\]
Then we have the representations (see \cite{Marchenco})

\begin{equation}\label{eq:2.1}
    \begin{cases}
        c(x, \lambda)=\cos{\lambda x} + \int\limits_{-x}^{x}K(x, t) \cos{\lambda t}dt, \\[3pt]
        s(x, \lambda)= \frac{\sin{\lambda x}}{\lambda} + \int\limits_{-x}^{x}K(x, t) \frac{\sin{\lambda t}}{\lambda}dt, 
    \end{cases}
\end{equation}
in which $K(x, t) \in C(\Omega) \cap W_1^1(\Omega)$, where 
\[\Omega=\{(x, t): 0 \leq x \leq 1, \; -x \leq t \leq x\},\] 
and $K(x, t)$ is the solution of the problem
\begin{equation}\label{eq:2.2}
    \begin{cases}
        K_{xx} - K_{tt} = q(x)K(x, t), \; \mbox{in} \;  \Omega\\[3pt]
        K(x, x) = \frac{1}{2}\int\limits_{0}^{x}q(t)dt, \: K(x, -x) = 0, \: x\in[0, 1]. 
    \end{cases}
\end{equation}
The following theorems are the main results of our previous work \cite{Biyarov} for the Sturm-Liouville operator 
\begin{theorem}[see Theorem 1.1. in \cite{Biyarov}]
\label{theorem:1}
The spectrum of $L_{DN}$ coincides with the spectrum of $L_{ND}$ (i.e. $\sigma(L_{DN})=\sigma(L_{ND})$) if and only if $q(x)=q(1-x)$ on $[0, 1]$.
\end{theorem} 

\begin{theorem}[see Theorem 1.2. in \cite{Biyarov}]
\label{theorem:2}
The spectrum of $L_D$ coincides with the spectrum of  $L_N$, except zero (i.e. $\sigma(L_{D})\setminus \{0\}=\sigma(L_{N}) \setminus \{0\}$), and $0\in \sigma(L_N)$ if and only if 
\[
q_1(x)=\Bigg(\int \limits_{1}^{x}q_2(t)dt\Bigg)^2, \tag{BB}\label{MyEqBB}
\]
 where  $q_1(x)=(q(x)+q(1-x))/2$ and  $q_2(x)=(q(x)-q(1-x))/2$ on $[0, 1]$.
\end{theorem}

Furthermore, we assume, without loss of generality, that $0\in\sigma(L_{N(\frac{1}{2})}).$
The following theorem is the main result of this paper
\begin{theorem}
\label{theorem:1.1}
The whole spectrum of  $L_P$, except the lowest, or the whole spectrum of $L_{AP}$  consist of the eigenvalues with geometric multiplicity two if and only if 
\[q(x)=q\Big(\frac{1}{2}-x\Big) \,\, \mbox{on} \,\,\,\, [0, \frac{1}{2}] \] 
or
\[q_1(x)=\Bigg(\int \limits_{\frac{1}{2}}^{x}q_2(t)dt\Bigg)^2, \tag{B}\label{MyEqB}\]
where $q_{1}(x)=(q(x)+q(\frac{1}{2}-x))/2$, $q_{2}(x)=(q(x)-q(\frac{1}{2}-x))/2$ on $\big[ 0, \frac{1}{2} \big]$.  
\end{theorem} 


\section{Proof of Theorem \ref{theorem:1.1}}
\label{sec:3_ProofTh1}

Let the whole spectrum of $L_{P}$, except the lowest, or the whole spectrum of $L_{AP}$ consist of the eigenvalues with geometric multiplicity two. Then all roots of 
\[ \Delta^2(\lambda)=4\]
be double roots, except lowest. Then Hill's equation has two linearly independent periodical or anti-periodical solutions of period $1$ or $2$, respectably, for all roots, except the lowest \cite[p.19]{Magnus}.
It is known that 
\[
\sigma(L_{P}) = \{ \lambda \in \mathbb{C}: \, s'(1, \lambda) + c(1, \lambda)= 2 \},
\]
\[
\sigma(L_{AP}) = \{ \lambda \in \mathbb{C}: \, s'(1, \lambda) + c(1, \lambda) = -2 \}.
\]
Then all of the eigenvalues  $\{ \lambda_{n} \}_0^{\infty} $ of $L_P$, except the lowest, are roots of the following system of equations
\[ 
    \begin{cases}
        \Delta(\lambda)= s'(1, \lambda) + c(1, \lambda) =2, \\[3pt]
        \Delta'(\lambda)=0,  
    \end{cases}
\]
and all of the eigenvalues $\{ \lambda'_{n} \}_1^{\infty} $ of $L_{AP}$ are roots of the following system of equations
\[ 
    \begin{cases}
        \Delta(\lambda)= s'(1, \lambda) + c(1, \lambda) =-2, \\[3pt]
        \Delta'(\lambda)=0.  
    \end{cases}
\]
For all $\{ \lambda_{n} \}_1^{\infty} $ we have two linearly independent eigenfunctions  $c(x, \lambda_n)$ and $s(x, \lambda_n)$  of $L_P$ with properties
\[c(1, \lambda_n) = s'(1, \lambda_n)=1, \;  c'(1, \lambda_n) = s(1, \lambda_n)=0, \]
and for all $\{ \lambda'_{n} \}_1^{\infty} $ we have two linearly independent eigenfunctions  $c(x, \lambda'_n)$ and $s(x, \lambda'_n)$  of $L_{AP}$ with properties
\[c(1, \lambda'_n) = s'(1, \lambda'_n)=-1, \; c'(1, \lambda'_n) = s(1, \lambda'_n)=0.\]
We construct the following pair of linearly independent solutions
\begin{equation}\label{eq:3.1}
    \begin{cases}
        y_1(x, \lambda)= s'(\frac{1}{2}, \lambda)c(x, \lambda) - c'(\frac{1}{2}, \lambda)s(x, \lambda), \\[6pt]
        y_2(x, \lambda)=  c(\frac{1}{2}, \lambda)s(x, \lambda) - s(\frac{1}{2}, \lambda)c(x, \lambda),  
    \end{cases}
\end{equation}
 with propoties 
\[y_1(\frac{1}{2}, \lambda) = y'_{2}(\frac{1}{2}, \lambda)=1, \; \; y'_{1}(\frac{1}{2}, \lambda) = y_{2}(\frac{1}{2}, \lambda)=0.\]
It is clear that for $\lambda=\lambda_n, \, n=1, 2, \dots$ these solutions will also be periodic of period 1, and for $\lambda=\lambdaэ'_n, \, n=1, 2, \dots $ they will be anti-periodic (of period 2).
Then the following two options are possible:
\\
(I). $c(x, \lambda_n)$ and $y_1(x, \lambda_n), \, n=1, 2, \dots$ are linearly dependent, i.e., Wronskian
\[W(c(x, \lambda_n), y_1(x, \lambda_n))=-c'\Big(\frac{1}{2}, \lambda_n\Big)=0,\] 
i.e., $\lambda_n\in\sigma(L_{N(\frac{1}{2})}), \, n=1, 2, \dots$
as well as $s(x, \lambda_n)$ and $y_2(x, \lambda_n), \, n=1, 2, \dots$ are linearly dependent, i.e., Wronskian
\[W(s(x, \lambda_n), y_2(x, \lambda_n))=s\Big(\frac{1}{2}, \lambda_n\Big)=0,\] 
i.e., $\lambda_n\in\sigma(L_{D(\frac{1}{2})}), \, n=1, 2, \dots$.
\\
(II).  $c(x, \lambda_n)$ and $y_2(x, \lambda_n), \, n=1, 2, \dots$ are linearly dependent, i.e., Wronskian
\[W(c(x, \lambda_n), y_2(x, \lambda_n))=c\Big(\frac{1}{2}, \lambda_n\Big)=0,\] 
i.e., $\lambda_n\in\sigma(L_{ND(\frac{1}{2})}), \, n=1, 2, \dots$,
as well as $s(x, \lambda_n)$ and $y_1(x, \lambda_n), \, n=1, 2, \dots$ are linearly dependent, i.e., 
\[W(s(x, \lambda_n), y_1(x, \lambda_n))=-s'\Big(\frac{1}{2}, \lambda_n\Big)=0,\] 
i.e., $\lambda_n\in\sigma(L_{DN(\frac{1}{2})}), \, n=1, 2, \dots$.

The case of $\lambda=\lambda'_n, \, n=1, 2, \dots$ gives the same result.
Then in case (I), by virtue of Theorem \ref{theorem:2} we have that the condition \eqref{MyEqB} holds on $[0, 1/2]$, and in case (II), by virtue of Theorem \ref{theorem:1}  we get that the condition $q(x)=q(1/2-x)$ holds on $[0, 1/2]$. Hence the necessary conditions were proven.

We shall now prove the sufficiency. Let condition \eqref{MyEqB} be satisfied  or $q(x)=q(1/2-x)$ on $[0, 1/2]$.
It is known that these conditions are satisfied simultaneously if and only if $q(x)\equiv0$. Let \eqref{MyEqB} is satisfied. Then by virtue of Theorem \ref{theorem:2}, we have 

\begin{equation}\label{eq:2.2}
    c'(\frac{1}{2}, \lambda) = - \lambda^{2}s(\frac{1}{2}, \lambda),\;\;\mbox{for all} \;\;\lambda\in \mathbb {C}.
\end{equation} 
Since $q(x)=q(x+1)$  on $[0, 1]$, we get that the $q(x+1)$ satisfies the codition  \eqref{MyEqB} as well. This follows from the fact that 
\[
\int\limits_{-\frac{1}{2}}^{\frac{1}{2}}\frac{q(t+1)-q(\frac{1}{2}-(t+1))}{2}dt=0.
\]
Then the following equality holds:
\begin{equation}\label{eq:2.3}
    c'(\frac{3}{2}, \lambda) = - \lambda^{2}s(\frac{3}{2}, \lambda)\;\;\mbox{for all} \;\; \lambda\in \mathbb {C}.
\end{equation} 
If $c(x, \lambda)$ and $s(x, \lambda)$ are the fundamental system of solutions to the equation \eqref{eq:1.1}, then $c(x+1, \lambda)$ and $s(x+1, \lambda)$ are also solutions of  \eqref{eq:1.1}. We find easily that
\[
 \begin{cases}
        c(x+1, \lambda)= c(1, \lambda)c(x, \lambda) + c'(1, \lambda)s(x, \lambda), \\
        s(x+1, \lambda)=  s(1, \lambda)c(x, \lambda) + s'(1, \lambda)s(x, \lambda).  
 \end{cases} 
\]
Then the following two cases will take place:

1).
 $c(x, \lambda)$ and $c(x+1, \lambda)$ are linearly dependent, i.e., Wronskian

\[W(c(x, \lambda), c(x+1, \lambda))=c'(1, \lambda)=0,\]
where we denote the zeros of this equation by  $\mu_n,\;n=0,1,2, \dots$.
As well as $s(x, \lambda)$ and $s(x+1, \lambda)$ are linearly dependent, i.e., Wronskian
\[W(s(x, \lambda), s(x+1, \lambda))=-s(1, \lambda)=0,\]
where $\lambda=\mu_n,\; n=1,2,\dots$.
Then for all  $\lambda=\mu_n,\; n=0,1,2,\dots$, we have the following

\begin{equation}\label{eq:2.4}
 \begin{cases}
        c(x+1, \lambda)= c(1, \lambda)c(x, \lambda), \\
        s(x+1, \lambda)=  s'(1, \lambda)s(x, \lambda).  
\end{cases} 
\end{equation} 

Let $s(\frac{1}{2}, \mu_n)\neq0$ is satisfied.  
From  \eqref{eq:2.4}, we have 
\[
\begin{cases}
        s(\frac{3}{2}, \lambda)= s'(1, \lambda)s(\frac{1}{2}, \lambda), \\
        c'(\frac{3}{2}, \lambda)=  c(1, \lambda)c'(\frac{1}{2}, \lambda).  
\end{cases} 
\]
By virtue   \eqref{eq:2.2} and  \eqref{eq:2.3}, we obtain that
\[
\begin{cases}
        s(\frac{3}{2}, \lambda)= s'(1, \lambda)s(\frac{1}{2}, \lambda), \\
        - \lambda^{2} s(\frac{3}{2}, \lambda)=- \lambda^{2} c(1, \lambda)s(\frac{1}{2}, \lambda).  
\end{cases} 
\]
Then, if $\lambda\neq0$, we have $s'(1, \lambda)= c(1, \lambda)$.

Let $s(\frac{1}{2}, \lambda)=0$ is satisfied.
From  \eqref{eq:2.4}, we get
\[
\frac{c(1,\mu_n)}{s'(1,\mu_n)}=\frac{c'(x+1,\mu_n)s(x,\mu_n)}{c'(x,\mu_n)s(x+1,\mu_n)}=\lim\limits_{\lambda\to\\\mu_n}\frac{c'(x+1,\lambda)s(x,\lambda)}{c'(x,\lambda)s(x+1,\lambda)}.\]
When $x$ tends to $\frac{1}{2}$, we obtain the following
\[
\frac{c(1,\mu_n)}{s'(1,\mu_n)}=\lim\limits_{x\to\\\frac{1}{2}}\Big(\lim\limits_{\lambda\to\\\mu_n}\frac{c'(x+1,\lambda)s(x,\lambda)}{c'(x,\lambda)s(x+1,\lambda)}\Big) \qquad\qquad\qquad\qquad\qquad\qquad\qquad
\]
\[
=\lim\limits_{\lambda\to\\\mu_n}\Big(\lim\limits_{x\to\\\frac{1}{2}}\frac{c'(x+1,\lambda)s(x,\lambda)}{c'(x,\lambda)s(x+1,\lambda)}\Big)=\lim\limits_{\lambda\to\\\mu_n}\frac{c'(\frac{3}{2},\lambda)s(\frac{1}{2},\lambda)}{c'(\frac{1}{2},\lambda)s(\frac{3}{2},\lambda)}.
\]
If $\lambda\neq0$, by virtue   \eqref{eq:2.2} and  \eqref{eq:2.3}, we get
$s'(1, \lambda)= c(1, \lambda)$.
Since 
\[s'(1, \lambda) c(1, \lambda)=1,\]
we have
\[s'(1, \lambda)= c(1, \lambda)=\pm1.\] 
Then we obtain the following system equations
\[ 
\begin{cases}
         \Delta^2(\lambda)=\big(s'(1, \lambda) + c(1, \lambda)\big)^2=4, \\[3pt]
        \Delta'(\lambda)=0,  
\end{cases}
\]
where $\lambda=\mu_n,\;n=1,2, \dots$, and the whole spectrum of $L_{P}$, except the lowest, or the whole spectrum of  $L_{AP}$ consist of the eigenvalues with geometric multiplicity two.

Let $q(x)=q(1/2-x)$ on $[0, 1/2]$. 
Then by virtue of Theorem \ref{theorem:1}, we get 
\begin{equation}\label{eq:2.5}
    c(\frac{1}{2}, \lambda) = s'(\frac{1}{2}, \lambda),\;\;\mbox{for all} \;\;\lambda\in \mathbb {C}.
\end{equation} 
Since $q(x)=q(x+1)$ on $[0, 1]$, we have

\begin{equation}\label{eq:2.6}
    c(\frac{3}{2}, \lambda) = s'(\frac{3}{2}, \lambda)\;\;\mbox{for all} \;\;\lambda\in \mathbb {C}.
\end{equation} 

Let $s'(\frac{1}{2}, \mu_n)\neq0$ is satisfied.
From  \eqref{eq:2.4}, we have 
\[
\begin{cases}
        c(\frac{3}{2}, \lambda)= c(1, \lambda)c(,\frac{1}{2} \lambda), \\
        s'(\frac{3}{2}, \lambda)=  s'(1, \lambda)s'(\frac{1}{2}, \lambda).  
    \end{cases}
    \] where $\lambda=\mu_n, \, n=0,1,2, \dots$.    
By virtue \eqref{eq:2.5} and \eqref{eq:2.6}, we get

\[
\begin{cases}
        s'(\frac{3}{2}, \lambda)= c(1, \lambda)s'(,\frac{1}{2} \lambda), \\
        s'(\frac{3}{2}, \lambda)=  s'(1, \lambda)s'(\frac{1}{2}, \lambda).  
    \end{cases}
    \]
Then we have that $s'(1, \lambda)= c(1, \lambda)$.

Let $s'(\frac{1}{2}, \lambda)=0$ is satisfied.
From  \eqref{eq:2.4}, we have

\[
\frac{c(1,\mu_n)}{s'(1,\mu_n)}=\frac{c(x+1,\mu_n)s'(x,\mu_n)}{c(x,\mu_n)s'(x+1,\mu_n)}=\lim\limits_{\lambda\to\\\mu_n}\frac{c(x+1,\lambda)s'(x,\lambda)}{c(x,\lambda)s'(x+1,\lambda)}.
\]
When $x$ tends to $\frac{1}{2}$, we get the following

\[
\frac{c(1,\mu_n)}{s'(1,\mu_n)}=\lim\limits_{x\to\\\frac{1}{2}}\Big(\lim\limits_{\lambda\to\\\mu_n}\frac{c(x+1,\lambda)s'(x,\lambda)}{c(x,\lambda)s'(x+1,\lambda)}\Big)\qquad\qquad\qquad\qquad\qquad\qquad\qquad
\]
\[=\lim\limits_{\lambda\to\\\mu_n}\Big(\lim\limits_{x\to\\\frac{1}{2}}\frac{c(x+1,\lambda)s'(x,\lambda)}{c(x,\lambda)s'(x+1,\lambda)}\Big)=\lim\limits_{\lambda\to\\\mu_n}\frac{c(\frac{3}{2},\lambda)s'(\frac{1}{2},\lambda)}{c(\frac{1}{2},\lambda)s'(\frac{3}{2},\lambda)}.
\]
By virtue   \eqref{eq:2.5} and  \eqref{eq:2.6}, we get
\[s'(1, \lambda)= c(1, \lambda).\]
Since 
\[s'(1, \lambda) c(1, \lambda)=1,\]
we have 
\[s'(1, \lambda)= c(1, \lambda)=\pm1.\]
Then we get the following system equations

\[ 
\begin{cases}
         \Delta^2(\lambda)=\big(s'(1, \lambda) + c(1, \lambda)\big)^2=4, \\[3pt]
        \Delta'(\lambda)=0,  
\end{cases}
\]
for all $\lambda=\mu_n,\;\; n=1,2, \dots$, and the whole spectrum of $L_{P}$, except the lowest, or the whole spectrum of $L_{AP}$ consist of the eigenvalues with geometric multiplicity two.

2).  $c(x, \lambda)$ and $s(x+1, \lambda)$ are linearly dependent, i.e., Wronskian

\[W(c(x, \lambda), s(x+1, \lambda))=s'(1, \lambda)=0,\]
where we denote the zeros of this equation by  $\mu'_n,\;n=1,2, \dots$. 
As well as $s(x, \lambda)$ and $c(x+1, \lambda)$ are linearly dependent, i.e., Wronskian
\[W(s(x, \lambda), s(x+1, \lambda))=-c(1, \lambda)=0,\]
where $\lambda=\mu'_n,\; n=1,2, \dots.$
Therefore,
\[\Delta(\lambda)=c(1,\lambda)+s'(1,\lambda)=0.\] 
Hence, in the case 2), there are no periodic or anti-periodic solutions. Thus, Theorem 1.3  is completely proved.

\begin{remark}\label{rem3.1}
In the particular case  when $q(x)=q(1-x)$ on $[0, 1]$, from Theorem \ref{theorem:1.1},  we have Theorem 5.1 and Theorem 5.2 in \cite{Biyarov}.
\end{remark}


\end{document}